\documentclass{amsart} 

\usepackage{latexsym,amssymb,amscd,eucal}

\newtheorem{thm}{Theorem}[section]
\newtheorem{lemma}[thm]{Lemma}
\newtheorem{prop}[thm]{Proposition}
\newtheorem{cor}[thm]{Corollary} 
\newtheorem{Def}[thm]{Definition} 
 
\newtheorem{con}[thm]{Conjecture}
\newtheorem{ax}{Axiom}

\numberwithin{equation}{section}

\newcommand{\BZ}{{\mathbb{Z}}}
\newcommand{\BQ}{{\mathbb{Q}}}
\newcommand{\BR}{{\mathbb{R}}}
\newcommand{\De}{{\mathbb{D}}}

\newcommand{\Si}{{\Sigma}}
\newcommand{\ra}{{\rightarrow}}
\newcommand{\la}{{\lambda}}
\newcommand{\el}{{\ell}}
\newcommand{\be}{{\beta}}
\DeclareMathOperator{\Sign}{Sign}
\DeclareMathOperator{\Trace}{Trace}
\DeclareMathOperator{\Mod}{Proj}
\DeclareMathOperator{\FMod}{Free}

\DeclareMathOperator{\id}{id}
\DeclareMathOperator{\Image}{Image}

\begin{document}

\title
{Integrality for TQFTs}
\author{Patrick M. Gilmer}
\address{Department of Mathematics\\
   Louisiana State University\\
    Baton Rouge, LA 70803\\
USA}
\email{gilmer@math.lsu.edu}
\thanks{partially supported by NSF-DMS-0203486}

\keywords{cobordism category, quantum invariant,  strong shift equivalence, cyclotomic}
\subjclass{57M}
\date{January 9, 2004}

\begin{abstract} We discuss  ways that the ring of coefficients for
a TQFT can be reduced if one restricts somewhat the allowed cobordisms.
When we apply these methods to a TQFT associated to $SO(3)$ at an odd
prime
$p,$ we obtain a functor from a somewhat restricted cobordism
category
to the category of free finitely generated modules over a ring of cyclotomic integers :
$\BZ [\zeta_{2p}], $
if $p \equiv -1 \pmod{4}$, and
$\BZ [\zeta_{4p}],$
if $p \equiv 1 \pmod{4},$
where $\zeta_k$ is a primitive $k$th root of unity.
We  study the quantum
invariants of  prime power order simple cyclic covers of 3-manifolds.  We  define new invariants arising from strong shift equivalence and integrality. Similar results are obtained for some other TQFTs  but the modules are only guaranteed to be  projective.
\end{abstract}

\maketitle
\section{Introduction}

A  2+1 dimensional TQFT  (topological quantum field theory) is a functor from a 2+1 dimensional cobordism category to the category of finitely generated projective modules over  some ring.
The cobordism category should have objects which  are closed surfaces (possibly empty) and   morphisms which are 3-dimensional cobordisms between closed surfaces.  Each cobordism has a source and a target surface. Here for the purposes of exposition, we ignore some technical details which will be dealt with in the main body of the paper.   As is traditional in this field, we use the  letter `V'  to denote  the functor on objects, and the letter `Z'  to denote the functor on  morphisms.
Thus we refer to a single functor by the pair $(V,Z)$.  If $\Si$ is a surface, then $V(\Si)$ is a module, and if $N$ is a cobordism between surfaces $Z(N)$ is a linear map between the respective modules.  The empty surface  is assigned  the ring itself.  In this way, a closed 3-manifold which may be viewed as morphism from the empty surface to itself is assigned  a scalar (that specifies the induced linear map between one dimensional free modules). This scalar is called  the quantum invariant (associated to the TQFT) of the closed  3-manifold.

If one rescales the quantum invariant so that a sphere takes the value one, then for many interesting TQFTs, the resulting quantum invariant of closed connected 3-manifolds takes values in  a ring of  cyclotomic integers.  In particular H. Murakami \cite{Mu1, Mu2}  showed this for homology spheres in  the  $SO(3)$  and $SU(2)$ theories at odd prime roots of unity.  This theorem was generalized by Masbaum and Roberts to include general closed connected manifolds containing colored links. Further integrality results of this type have been obtained by Masbaum and Wenzl \cite{MW}, as well as Takata and Yokota \cite{TY}.  Le has obtained very general integrality results of this kind \cite{Le}.
 We will show that this integrality of quantum invariants implies that much of the machinery of the TQFT exists over the corresponding ring of integers.

We will call a  cobordism in the above category ``targeted''  if each component of the cobordism meets
at least one  component of the target surface. The targeted cobordisms are closed under composition and so form a subcategory.  If the ring of coefficients is a  subring of a number field, we may extend the coefficients to the number field.  We also suppose that the rescaled quantum invariants of closed manifolds  lie in  the ring of  integers of this number field. In this situation, we assign a finitely generated projective module over the ring of integers to each surface.  

A targeted morphism from the empty set to a given surface induces via the  TQFT
over the number field,  a map from the field to the vector space associated to the surface.
The image of one in the field is a vector in this vector space.   We call this vector a targeted vacuum state. The new module  that  we associate to a surface is the module spanned,  over the ring of integers of the field,  by all targeted vacuum states.  In this way,  the original TQFT functor restricted to the targeted subcategory factors through a functor which takes values in the category of finitely generated projective modules over the ring of integers.

This procedure, which relies on the fact that a ring of integers in a number field is a Dedekind domain, is described  in \S2 in a more general context. We consider the problem of reducing the ring of coefficients from a commutative integral domain with unit to a subring which is Dedekind.  We don't actually require that our initial functor satisfy the tensor product axiom. As there are many  different axiomatic approaches to TQFT, we spell out exactly which properties are needed for our procedure to go through. We believe this clarifies the argument. In particular we consider more general cobordism categories than the 3-dimensional one discussed above. In fact we make use of Turaev's notions of $\mathfrak {A}$ and $\mathfrak {B} $ space structures.
In this axiomatic context, we prove the integrality of  the (unrescaled) quantum invariants of 3-manifolds with first Betti number nonzero. We obtain restrictions on the quantum invariants of 3-manifolds with simple cyclic group actions of prime order. In this general context, we also  obtain the integrality of Turaev-Viro polynomials of knots \cite{1dim, sat}.  

The  modules obtained in \S 2 are not necessarily free. We will 
 show, for a version of the $SO(3)$ theories associated to an odd prime $p$,   that
the modules over cyclotomic rings of integers obtained in this way are free.  To show that the above-mentioned modules over cyclotomic rings of integers are free,  we  use Lemma \ref{inj}.  It says, assuming $p$ is a prime, that if a finitely generated projective  $\BZ [\zeta_{2p} ]$-module becomes free after we localize by inverting $p$, then it was free already.  

We  begin with the TQFTs  of Blanchet, Habegger, Masbaum, and Vogel   \cite{BHMV2}  which depend on a choice of an integer $p\ge 3.$  When $p$ is odd, this is a version of the TQFT associated to  $SO(3).$   Moreover the modules associated to a surface are  free in these theories.
 We adapt these theories by replacing the cobordism category used  in \cite{BHMV2} where objects and morphisms possess the extra data of a  $p_1$-structure with another cobordism category, $\mathcal{C}$, where the surfaces possess the extra structure of a Lagrangian subspace of the rational first homology, and the cobordisms possess the extra structure of an  integral weight.  \S 3 describes this subcategory. The idea  to use this kind of data to resolve the framing anomaly  in TQFT is due to  Walker \cite{Wa}. It was further developed by  Turaev \cite{Tu}.  
In \S 4, we  describe the resulting TQFTSs  which will be denoted  $(V_p,Z_p)$. These theories may be defined over $\BZ [\zeta_{2p},\frac 1 p ]$ when $p \equiv 0 \text{ or }  -1  \pmod{4}, $   and over $\BZ [\zeta_{4p},\frac 1 p ]$ when $p \equiv 1 \text{ or } 2 \pmod{4}$. We let $K_p$ denote these rings, and $\Bbb{D}_p$  their rings of integers. 

In \S 5, we apply the construction of \S 2 to the $(V_p,Z_p)$ theory, in the case where
$p$ is an odd prime or twice an odd prime, obtaining functors $(S_p,Z_p)$ to the category of finitely generated projective modules over $\Bbb{D}_p$.  In \S6, $S_p(\Si)$ is then shown to be free when $p$ is prime,  and $p \equiv -1 \pmod{4}$, using Lemma \ref{inj} .

In \S7 we describe an `index two' subcategory of the cobordism category. 
It is defined by  constraining the weights of cobordisms by a congruence modulo two.
In \S 8, we use this subcategory to  reduce the coefficients from $\BZ [\zeta_{4p},\frac 1 p ]$ to  $\BZ [\zeta_{2p},\frac 1 p ]$ coefficients in the case  $p \equiv 1 \text{ or } 2 \pmod{4}$.  The theories $(V_p,Z_p)$, when restricted to the cobordism subcategory of \S 6,  yield functors $(V^+_p,Z_p)$ over $\BZ [\zeta_{2p},\frac 1 p ].$ 

 In \S 9, we apply the method of \S 2 to $(V^+_p,Z_p)$
obtaining functors $(S^+_p,Z_p),$ when $p$ is either an odd prime  and $p \equiv 1 \pmod{4}$ or when $p$ is twice an odd prime. $S^+_p(\Si)$ is a finitely generated projective module over
$\Bbb{D}_p$.
We then use Lemma \ref{inj}  to prove that
$S_p^+(\Si)$ is free in the case $p$ is prime, and
$p \equiv 1 \pmod{4}.$ The freeness of
$S_p(\Si)$  in the case $p \equiv 1 \pmod{4}$ follows from that of $S_p^+(\Si).$

The main goal of  this paper  is the following theorem which follows immediately from Theorems \ref{main} 6.2, and 9.1. 
Let  $\mathcal{C'}$ denote the category with the same objects as   $\mathcal{C}$ but with only targeted morphisms. Let $\FMod(K)$ denote  the category of 
free finitely generated $K$-modules.
\begin{thm}\label{simple} Let $p$ be an odd prime.  $(S_p,Z_p)$ is a functor from $\mathcal{C'}$  to $\FMod(\De_p)$. The composition of   $(S_p,Z_p)$ with the localization functor from $\FMod(\De_p)$ to $\FMod(K_p)$  is naturally isomorphic to the restriction of $(V_p,Z_p)$ to $\mathcal{C'}.$
\end{thm}

It is very desirable to have explicit bases for the free modules $S_p(\Sigma)$ for each odd prime $p.$ In joint work with  Masbaum, and  van Wamelen \cite{GMvW},  we have found such bases for surfaces of genus one and two without any colored banded points.  This has just recently been extended by Masbaum and myself to all surfaces without colored banded points   \cite{GM}.   Kerler  has announced an integral version of the SO(3)  theory at p=5, with explicit bases for surfaces  of any genus without colored banded points.

In \S 10 we discuss applications to strong shift equivalence class invariants of infinite cyclic covers. In \S 11 we study applications to the $SO(3)$ quantum invariants of 3-manifolds with simple prime power order cyclic covers. Theorem \ref{simplecong} was used in joint work \cite{GKP} with Kania-Bartoszynska and  Przytycki  to obtain restrictions on the $SO(3)$ quantum invariants  of  periodic homology spheres. Chen and  Le \cite{CL}  have generalized \cite{GKP}  to projective quantum invariants associated to simple complex Lie groups. This generalization uses some results of section \S 2  which were given in an axiomatic context, precisely with this kind of application in mind.  We close the paper with a conjecture about the divisibility of quantum invariants.

Theorems \ref{main} 6.2, 7.2, and 9.1 were first announced at Knots 2000 in August 2000. At that time, Theorems 6.2 and  9.1  were known only for $p < 23,$ when the class number of $\BZ[\zeta_p]$ is one. Theorems 6.2 and  9.1 without this restriction were announced at the ``Low dimensional topology'' special session at the A.M.S. meeting in New Orleans in January 2001.  Conjecture 12.1 was made there as well.  Theorem 2.5 is stated here with a slightly different restricted cobordism category $\mathcal{C}'$ and a slightly different module $S(\Sigma)$ associated to a disconnected surface $\Sigma.$  One reason we made these modifications is so that Proposition 2.14 would hold.

We would like to thank Stavros Garoufalidis, Thomas Kerler, Gregor Masbaum, Jorge Morales, Khaled Qazaqzeh, and
Paul van Wamelen for helpful discussions.

We use  $\beta_i(X,Y)$ to denote the ith Betti number of a pair $(X,Y).$

\section{Reduction of the ring of coefficients for targeted morphisms}

In this section we work axiomatically. We will use Turaev's concepts of
 space structure, and  of cobordism theory
\cite[III]{Tu}.
Let $\mathfrak {A}$ and $\mathfrak {B} $ be space structures compatible with
disjoint union, and let  $\mathfrak {A}$ be involutive. Suppose
$(\mathfrak {B}, \mathfrak {A})$ forms a cobordism theory.

We define a category $\mathcal{C}$ whose objects are $\mathfrak {A}$-spaces. A
morphism of $\mathcal{C}$ from $\Sigma_1$
to $\Sigma_2$ is an equivalence class of a pair consisting of a
$\mathfrak {B}$-space $M,$
and  an $\mathfrak {A}$-homeomorphism (called a boundary identification)
of $\partial M$ to
$-\Sigma_1
\sqcup \Sigma_2.$  We let $\partial_{-} M$ denote the part  of $\partial M$ that  is identified with
$\Sigma_1,$ and $\partial_{+} M$ denote the part  of $\partial M$ that is identified to
$\Sigma_2.$
Two $\mathfrak {B}$-spaces with boundary
identification are considered equivalent if there is a $\mathfrak
{B}$-homeomorphism between them that preserves the boundary
identification. Composition is defined by gluing of $\mathfrak
{B}$-spaces. A category formed in this way is called a cobordism
category.

Let $K$ be a commutative integral domain with identity. Let $\Mod(K)$ denote
the category of finitely generated projective $K$-modules.

Let $(V,Z)$ denote a functor that assigns to objects of $\mathcal{C}$ objects of $\Mod(K)$ and assigns
to morphisms of $\mathcal{C}$ morphisms of $\Mod(K)$.
$$\text{ a object of }  \mathcal{C} \mapsto V(\text{object}), \text{a finitely generated projective}  \ 
K \text {-module} $$
$$\text{ a morphism of } \mathcal{C} \mapsto Z(\text{morphism})
	\text{,  a  } K \text{-linear map}
 $$
{\em We will assume that $(V,Z)$
satisfies the following three axioms.}  We note that the anomaly-free  TQFT
that Turaev constructs from a modular category will satisfies these axioms as well as have the trace property, given in  Definition \ref{trace}.
The same holds for a 3-dimensional TQFT in sense of Blanchet-Habegger-Masbaum-Vogel that satisfies their surgery axioms (S0) and (S1).

\begin{ax}  \label{empty} $V(\emptyset)= K.$ \end{ax}

If $ N:\emptyset \rightarrow \Sigma,$ then
$Z(N)(1) \in V(\Sigma).$  This element is called a $\mathcal{C}$-vacuum state and
is denoted by  $[N].$ If $N$ is connected, we say $[N]$ is a connected
$\mathcal{C}$-vacuum state for $\Sigma.$ Functors satisfying the following axiom without the connectivity condition are called
`cobordism generated' by \cite{BHMV2}, and `nondegenerate' by \cite{Tu}.

\begin{ax}  \label{cg} For each object $\Sigma,$  $V(\Sigma)$ is finitely generated by
connected $\mathcal{C}$-vacuum states.\end{ax}

If $M:\emptyset \rightarrow \emptyset,$ then  $[M] \in K.$  Instead of writing $M:\emptyset \rightarrow \emptyset,$ we will sometimes simply write `$M$'   is closed.   When $M$ is closed, we will  follow the traditional notation and write $\langle M \rangle$ for $[M].$

\begin{ax}  \label{bform}
There is a bilinear form
 $$v_{\Sigma}:V(\Sigma ) \times V(-\Sigma) \rightarrow K,
 \quad \text{satisfying} \quad v_{\Sigma}([N_{1}],[N_{2}])= \langle N_{1} \cup_{\Sigma}N_{2} \rangle$$
$\text{where} \quad N_{1}: \emptyset \rightarrow \Sigma \quad \text{and} \quad
N_{2}: \emptyset \rightarrow -\Sigma.$
Moreover this form has an injective
 adjoint homomorphism:
 $$V(\Sigma ) \rightarrow
\hom (V(-\Sigma),K) $$ \end{ax}

\begin{Def} We say $N: \Sigma \rightarrow
\Sigma'$  is targeted if $\beta _0(N,\Sigma')=0, $ i.e. every component of $N$ includes some component of 
$\Sigma'.$ \end{Def}

Note the if $N: \Sigma \rightarrow
\Sigma'$  is targeted  and $N$ is nonempty,  then
$\Sigma'$ is nonempty. Also $N$ cannot be closed if it is targeted.
 Note that the composition of  two targeted morphisms is  again targeted.

\begin{Def}Let $\mathcal{C}'$ be the subcategory of $\mathcal{C},$ with the same objects as $\mathcal{C},$ but with only targeted morphisms. \end{Def}

We say $[N]$ is a $\mathcal{C}'$-vacuum state if $N:\emptyset \ra \Si$ is a morphism of $\mathcal{C}'.$

{ \em We now assume that $K$ contains a Dedekind domain $ {\De}.$ }

\begin{Def}
We say $(V,Z)$ is almost $ {\De}$-integral if there exists some 
 $\mathcal{D}\in K$ such that
$\mathcal{D} \langle M \rangle \in   {\De},$ for every  connected closed morphism $M.$
\end{Def}

\begin{Def}  If $\Sigma \ne \emptyset, $  we define
$S(\Sigma)$ to be the $ {\De}$-submodule of $V(\Sigma)$ generated
by all $\mathcal{C}'$-vacuum states for $\Sigma.$   We define
 $S(\emptyset)= {\De}.$
\end{Def}

\begin{thm}\label{main}
  If $(V,Z)$ is almost $ {\De}$-integral, then

 a) $S(\Sigma)$ is a finitely generated projective $ {\De}$-module.

 b) $S(\Sigma) \otimes_{ {\De} } K$ is naturally isomorphic to  $V(\Sigma)$ as $K$-modules

 c) If $N: \Sigma \rightarrow
\Sigma'$  is targeted then
 $Z(N):S(\Sigma) \rightarrow S(\Sigma').$

In summary, we have a commutative (up to natural isomorphism) diagram of categories and
functors:
\[
\begin{CD}
\mathcal{C}' @>>> \mathcal{C} \\
 @V(S,Z)VV                       @V(\text{$V$},Z)VV \\
\Mod( {\De}) @>\otimes_{K}>> \Mod(K)
\end{CD}
\]
\end{thm}

\begin{proof} By Axiom 2, we can
  pick a finite set of generators:
  $\{ [G_{i}] \}_{1 \le i\le n}$ for $V(-\Sigma)$  where
 $G_{i}: \emptyset  \rightarrow -\Sigma$ are connected.
 Let $T$ denote the collection of all targeted morphisms from the empty set to $\Si$.  
For any ring $R$,  let $R^{(T)}$ denote the $R$-module of
formal linear combinations of elements of $T$ with coefficients in $R$, and  with only 
finitely many non-zero coefficients.  Define $\psi: \De^{(T)} \ra \De^n$  on generators by
$$\psi(N) = \ \mathcal{D}  ( v_{\Sigma}([N],[G_{1}]), \cdots,
 v_{\Sigma}([N],[G_{n}])).$$  
 $\psi$ takes values in $\De^n$, by the 
almost ${\De}$-integrality property. We also have a surjective $\De$-module map
from $ \phi :\De^{(T)} \ra S(\Sigma)$ which sends $N$ to  $[N].$  $\psi$ factors through  $\phi$, and induces a surjection between  $S(\Sigma)$ and the image of  $\psi$.  This surjection is injective by Axiom 3. Thus we have an isomorphism between $S(\Si)$ and the image of $\psi.$ 
But a submodule of ${\De}^n$ is necessarily finitely generated
and projective \cite[Prop10.12]{Ja}. This proves statement a).

Let $U_{\De}$ be the kernel of $\phi$ and  let $U_{K}$ be the kernel of the canonical map 
$K^{(T)} \ra V(\Sigma)$ which sends $N$ to $[N]$.   Then we have  short exact sequences:
  \begin{equation} \label{exactS} 0 \ra U_{\De}  \ra \De^{(T)} \ra S(\Si) \ra 0 . \end{equation}
  \begin{equation} \label{exactV} 0 \ra U_{K}    \ra K^{(T)}  \ra V(\Si) \ra 0  .\end{equation}
 
 We will show that  \begin{equation}\label{claim} U_{\De} \otimes K = U_{K} \end{equation}
  As $S(\Si)$ is projective,  the sequence \eqref{exactS} splits \cite[Prop 3.10]{Ja}, and we may obtain an exact sequence by tensoring  equation \eqref{exactS}
 with $K$:
\begin{equation} \label{exactSK}
 0 \ra U_{\De} \otimes K  \ra \De^{(T)}  \otimes K \ra S(\Si)\otimes K  \ra 0 .\end{equation}
 Comparing the exact sequences \eqref{exactV} and \eqref{exactSK} and using equation  \eqref{claim},  we obtain the natural isomorphism 
 of claim b).   
 We now show that equation \eqref{claim} holds. By the proof of statement a) above, $U_{\De}$  is  the kernel of $\psi$.
Similarly, $U_{K}$  is isomorphic to the kernel of  $\psi \otimes \id_K:  K^{(T)}  \ra K^n.$ So we have exact sequences

\begin{equation}  \label{exactimD} 0 \ra U_{\De}  \ra \De^{(T)} \ra \Image(\psi) \ra 0 . \end{equation}
\begin{equation} \label{exactimK} 0 \ra U_{K}   \ra K^{(T)}  \ra  \Image(\psi \otimes \id_K)= \Image(\psi) \otimes K \ra 0  .\end{equation}
 As  $\Image(\psi)$ is projective,  we may, as above, obtain an exact sequence by tensoring  
 sequence \eqref{exactimD}
 with $K$:
  \begin{equation} \label{exactimDK}
 0 \ra U_{\De} \otimes K  \ra \De^{(T)}  \otimes K \ra \Image(\psi)  \otimes K \ra 0 .\end{equation}
 
 Comparing the exact sequences \eqref{exactimK}  and \eqref{exactimDK}, 
 we see that $U_{\De} \otimes K = U_{K} $. 
Statement c) follows as the composition of two targeted morphisms is again targeted.
\end{proof}

\begin{cor}
If $(V,Z)$ is almost $ {\De}$-integral, and $ {\De}$ is a PID, then $S(\Sigma)\approx  {\De}^k$ for some k.
\end{cor}

 The following corollary was suggested by \cite{rep}.

\begin{cor} If $(V,Z)$ is almost $ {\De}$-integral  and
 $V(\Sigma)$ is free, then there is an ideal $I$ of $ {\De}$ and
  a basis of $V(\Sigma)$ such that for every targeted morphism $N:
\Sigma \rightarrow
\Sigma,$   $Z(N)$ is given by a matrix with
entries in the fractional ideal $I^{-1} = \{r \in K |r I \subset  {\De} \}.$
\end{cor}

\begin{proof} By \cite[Theorem 10.14]{Ja},
$S(\Sigma)\approx  {\De}^k \oplus I,$ where $I$ is some ideal in $ {\De}.$
 Thus $V(\Sigma)\approx K^{k+1}.$ It is easy to check that the
corresponding basis for
$V(\Sigma)$ has the claimed property.
\end{proof}

We define the characteristic polynomial of an endomorphism 
${ \mathcal Z}: { \mathcal S} \ra { \mathcal S}$
in $ \Mod( {\De})$  to be the characteristic polynomial of the endomorphism 
${ \mathcal Z}: { \mathcal S } \otimes {\mathcal K}  \ra { \mathcal S} \otimes {\mathcal K} ,$ 
  where ${\mathcal{K}}$ is the field of fractions of $ {\De}$. The following corollary together  
  with \S 7 and \S8 answers \cite[Conjecture 1 p.256]{1dim} when $p$ is an odd prime or twice an odd prime.

\begin{cor}
 \label{char}
 If $(V,Z)$ is almost $ {\De}$-integral, and
   $N: \Sigma \rightarrow
\Sigma$ is targeted, then the characteristic polynomial for $Z(N)$
lies in $ {\De}[x].$
\end{cor}
\begin{proof}
There is an ideal $J$ in $ {\De}$ such that $S(\Sigma) \oplus J \approx
 {\De}^l$
\cite[Theorem 10.14]{Ja}.
So the characteristic polynomial of $Z(N)\oplus \text{0}_J$ lies in
$ {\De}[x].$
The characteristic polynomial of $Z(N)\oplus \text{0}_J$ is $x$
times the characteristic polynomial of $Z(N).$
\end{proof}

\begin{Def} Suppose $N: \Sigma \rightarrow \Si$ is a morphism in $\mathcal{C}$.
Let $C(N),$ the closure of $N,$ be the closed $\mathfrak B$-space obtained from $N$ by gluing
$\partial_-N$ to $\partial_+N.$ \end{Def}

\begin{Def}  \label{trace}  $(V,Z)$ is said to have the trace property if for any $N: \Sigma \rightarrow \Si,$ $<C(N)>= {\Trace}(Z(N)).$
\end{Def}

Suppose $M$ is a closed connected morphism, and $\pi:M_{\infty} \rightarrow M$ is a connected infinite cyclic covering space of $M.$ We would like to be able to choose a bicollared $\mathfrak {A}$-space $\Si$ (called a slit)  in $M$  such that :
\begin{itemize}

\item
$M \setminus \Si $ is the interior of a $\mathfrak {B}$-space $N$ with $\partial N= X_1 \sqcup -X_2$ where each $X_i$ is equipped with a $\mathfrak {A}$-homeomorphism to $\Si.$

\item There is a section $s$ of $\pi_{|\pi^{-1}(M \setminus \Si) }$

\item $M$ is the closure of $N$ when regarded as an endomorphism of $\Si.$
\end{itemize}

\begin{Def}  $\mathcal{C}$ is said to have the slitting property if we can choose a slit $\Si$ for any infinite cyclic cover of a
closed connected $\mathfrak {B}$-space $M$. We say $\mathcal{C}$  has the connected slitting property, if
 we can choose $\Si,$ so that, in addition, $N$ is connected.\end{Def}

From now on, suppose that every $\mathfrak {B}$-space is locally path connected and semi-locally 1-connected. We do this so that one can form covering spaces corresponding to subgroups of the fundamental group.

\begin{cor} \label{betti} Suppose $(V,Z)$ is almost $ {\De}$-integral, and has the trace property. If $\mathcal{C}$ has the connected slitting property, and
 $M$ is a closed connected $\mathfrak B$-space  with first Betti number nonzero, then $<M>
\in  {\De}.$
\end{cor}
\begin{proof}   $M$  is the closure of an endomorphism in $\mathcal{C}'.$ So  $ <M>=
 {\Trace}\  Z(N).$ This is a coefficient of the characteristic polynomial, and so, by
Corollary \ref{char}, lies in $ {\De}.$
\end{proof}

A $\mathbb {Z}_d$-covering space $\tilde X$ of $X$ is said to be simple if it is classified by a map $H_1(X) \rightarrow \mathbb{ Z}_d$ that factors through an epimorphism $\chi$ to $\mathbb{ Z}.$ If $M$ is closed,  $\tilde M$ is a  simple $\mathbb{ Z}_d$-covering space of $M,$ and $\mathcal{C}$ has the  slitting property, then $\tilde M$
inherits an $\mathfrak A$-structure once we make a choice of slit for $\chi.$ We will say this is a compatible $\mathfrak A$-structure. (In  the cases of interest, this is independent of the choice of slit.)

Note that if  $M$ is closed, and  $\tilde M$ is a  simple $\mathbb{ Z}_d$-covering space of $M,$ then both $M$ and $\tilde M$ have first Betti number greater than zero and so $<M>$ and $<\tilde M>$ already lie in $ {\De}$ by Corollary \ref{betti}.   

\begin{cor} \label{scover}  Suppose $ {\De}$ is a ring of algebraic integers, $(V,Z)$ is almost $ {\De}$-integral, and has the trace property, and
$\mathcal{C}$ has the connected slitting property.   Let $M$ be closed. Let $\tilde M$ be a  simple $\mathbb{ Z}_{d}$-covering space of $M$ with compatible $\mathfrak A$-structure. Suppose $d=r^s,$ where $r$ is a prime and $ {\De}$ is a ring of algebraic integers. We have:
$$ <\tilde M>- < M>^{d} \in r {\De}.$$
\end{cor}

\begin{proof}
We have $N:\Si \rightarrow \Si,$ so that $M=C(N),$ and $\tilde
M=C(N^d).$ Thus $Z(N): S(\Si) \ra S(\Si),$  $<M>=
  {\Trace}(Z(N)),$ and $<\tilde M>= \text{Trace}(Z(N^d)) =
 {\Trace}(Z(N)^d).$ Choose $J$ as in the proof of Corollary
\ref{char}. Let $\mathcal{Z}=Z(N)\oplus \text{0}_J.$ $\mathcal{Z}$ is
then an endomorphism of a free $ {\De}$-module. Moreover $<M>=
{\Trace}(\mathcal{Z}),$ and $<\tilde M>= {\Trace}(\mathcal{
Z}^d).$ But ${\Trace}(\mathcal{Z})= \sum_i \la_i,$   where
$\{\la_i\}$ denotes the list the eigenvalues of $\mathcal{Z},$ with
multiplicity. Each $\la_i$ is itself an algebraic integer as it
satisfies a polynomial in $ {\De}[x],$ namely the characteristic
polynomial of $\mathcal{Z}.$ Moreover ${\Trace}(\mathcal{Z}^d)=
\sum_i \la_i^d.$ The result follows as $\sum_i \la_i^d-(\sum_i
\la_i)^d$ is $r$ times an algebraic integer.
\end{proof}

We close this section with some further general properties of the  resulting functor
$(S,Z).$

\begin{prop} \label{tensor} There is a natural map:
$$S(\Si_1) \otimes S(\Si_2) \ra S(\Si_1 \sqcup \Si_2).$$ \end{prop}

\begin{prop} One has a bilinear pairing $s_{\Sigma}:S(\Sigma) \times
S(-\Sigma) \rightarrow  {\De}$ with an injective adjoint $S(\Sigma)
\ra \hom(S(-\Sigma), {\De}) $
satisfying $
 s_{\Sigma}([N_{1}],[N_{2}])= \mathcal{D}^{\beta_0(\Si)}  \langle N_{1}
\cup_{\Sigma}N_{2} \rangle.$ Suppose  $(\mathfrak {B}, \mathfrak {A})$ is an involutive cobordism theory \cite[III.2.8.3]{Tu} and that
$K$ is equipped with an involution  that
sends $<-M>$ to $<M> $ for all closed morphisms $M$, fixes $ \mathcal{D}$ and preserves
 $ {\De}$  then
 there is  a  Hermitian form
$ ( \  , \  )_{\Sigma}: S(\Sigma) \otimes_ {\De} S(\Sigma) \rightarrow  {\De}$ given by $ ( [N_1],[N_2]  )_{\Sigma}= \mathcal{D}^{\beta_0(\Si)}  \langle N_1 \cup_{\Sigma} -N_2 \rangle$.
 \end{prop}

 \section{some 3-dimensional involutive cobordism theories}

We  recall the space structures $\mathfrak A^e,$
${\mathfrak B}^w,$ and the resulting  3-dimensional involutive cobordism theory
$({\mathfrak B}^w, {\mathfrak A}^e)$ as described by Turaev \cite{Tu}. The idea
to use Lagrangians and integral weights as extra data to resolve anomalies in TQFT is due to Walker \cite{Wa}. We are influenced by both sources.

To simplify the exposition, we work initially with bare manifolds which are not decorated with colored banded trivalent graphs. These will be added later.

A ${\mathfrak A}^e$-space or an extended surface \cite[IV.6.1]{Tu} is an oriented compact surface $\Sigma$
together with a choice of a Lagrangian subspace $\el({\Sigma})$
for $H_1(\Sigma,\BQ )$ with respect to the intersection pairing. i.e., $\el({\Sigma})$ is a subspace of $H_1(\Sigma,\BQ)$
that is self-annihilating with respect to the intersection
pairing on  $H_1(\Sigma,\BQ).$ An e-homeomorphism between
extended surfaces is a homeomorphism that preserves the orientation and the Lagrangian subspaces.  Note that we have departed from Turaev's treatment by
taking Lagrangian subspaces of $H_1(\Sigma,\BQ )$ rather than
$H_1(\Sigma,\BR ).$ This is necessary for our description of a cobordism subcategory in \S4.

As described in \cite[Section 1]{Wa},
one should think of $\el({\Sigma})$ as a remnant of a 3-manifold $P$
whose boundary is $-\Si.$  The kernel of the map on
$H_1(\Si, \BQ) \ra H_1(P, \BQ)$ induced by the inclusion  is a Lagrangian subspace of $H_1(\Sigma,\BQ ).$ We will call such a $P$ a   $\Si$-coboundary.

A ${\mathfrak B}^w$-space or a weighted extended 3-manifold \cite[IV.9.1]{Tu} is an oriented compact 3-manifold $N$ equipped with
an integer weight $w(N),$  whose boundary is equipped with an
${\mathfrak A}^e$-structure. The empty 3-manifold $\emptyset$ is only allowed the weight zero.

As described in \cite[Section 1]{Wa},
one should think of the weight
of a  $\mathfrak B$-space $N$ as the signature of some background
4-manifold whose boundary is the result of capping $N$ off with some
$\partial N$-coboundary.

Suppose $N$ is a weighted, extended 3-manifold and
$\partial N= X \sqcup Y  \sqcup Z,$ and $g:X \ra -Y$ is a e-homeomorphism.
Here $X,$ $Y,$ and $Z$ are extended surfaces in their own rights. In this situation,  the weight of  the extended 3-manifold resulting from gluing $X$ to
$Y$ by $g$ is defined \cite{Wa} \cite{Tu}.    Its weight  is defined using $w(N),$
$\el(X),$ $\el(Y)$ $\el(Z)$ , $g_*$,  and the homological properties of  $N.$ The formula involves the  Maslov index  (or Wall non-additivity correction term) of three Lagrangian subspaces. 
Actually the formulas \cite[1.12]{Wa}  and  \cite[p 214]{Tu} differ by the sign in front of the Maslov index term.\footnote{In comparing these formulas, we must note that
Turaev's  Maslov index  is minus  Wall's correction factor, which is used by Walker, but we also have to notice that the order of the three  Lagrangians differs by an odd permutation in the two formulas.}  We follow Walker's definition, as then the following signature interpretation of the weight, due  to Walker,  will hold.  
 Let $P_X$ be an
$X$-coboundary,
 and $P_Z$ be a $Z$-coboundary. Pick $P_Y,$  a  $Y$-coboundary such that $g$ extends to an orientation preserving homeomorphism of $P_X$ to $-P_Y.$
 Let $W$ be an oriented 4-manifold with boundary
$N \cup_{X} P_X \cup_{Y} P_Y \cup_{Z} P_Z$ whose signature is $w(N).$
Let $W'$ be the result of gluing $P_X$ to $P_Y.$ Then $w(N')$ is the signature of $W'.$
We use  the signature interpretation  to prove Proposition \ref{defect},  Theorem \ref{even}, and Proposition \ref{evenslice}.

The sign convention chosen above has the effect of making the surgery formula for the quantum invariant of a closed weighted manifold especially simple when the signature of the linking matrix for the surgery description is exactly the weight. There are equivalences between  the cobordism categories resulting from either  choice of sign convention. These equivalences are defined by multiplying weights by minus one.

 Recall the Hirzebruch signature defect \cite{hirz} which has a somewhat simpler definition in the case of simple covers:
If $\tilde M \rightarrow M$ is a simple $\BZ_d$ covering of closed 3-manifolds, it is the boundary of a simple  $\BZ_d$ covering space $\tilde W \rightarrow W$ of $4$-manifolds.  One  defines $\text{def}( \tilde M \ra M) = d \  \Sign( W) -\Sign (\tilde W).$

The following Proposition follows directly from the definition and the 4-manifold interpretation of induced weights. An analogous result for any finite regular cover with the $\sigma$-invariant of $p_1$-structures playing the role played here by weights of extended weighted 3-manifolds is given in \cite[Lemma 4]{per}.

\begin{prop}
 \label{defect}
If $\tilde M \rightarrow M$ is a simple $\BZ_{d}$ covering of weighted 3-manifolds then
$w(\tilde M)= d w(M) -\text{ def }(\tilde M \rightarrow M).$
\end{prop}

From now on, we simply let $\mathcal{C}$ denote the cobordism category associated to
the cobordism theory $({\mathfrak A}^e,{\mathfrak B}^w)$ as in \S 2. Note $\text{Id}_{\Si}:
\Si \ra \Si$ is given by $\Si \times I$ weighted zero.

We now consider some associated decorated cobordism theories.
Let $p$ be an integer greater than or equal to $3.$   A
${\mathfrak A}^e_p$-space
is a ${\mathfrak A}^e$ space with a $p$-colored banded set of points. Here $p$-colored means labelled with nonnegative integers less than $p-1,$ if $p$ is odd, and less than $\frac p 2 -1$ if $p$ is even.  
A ${\mathfrak B}^w_{p}$-space is a ${\mathfrak B}^w$ space with some 
$p$-colored  banded trivalent graph \cite[\S 4]{BHMV2} .  
Let $\mathcal{C}_p$ denote the cobordism category associated to the cobordism theory
$( {\mathfrak A}^e_p, {\mathfrak B}^w_p ).$  Let  ${\mathcal{C}_{\breve p}}$ denote the subcategory of $\mathcal{C}_p$ where all objects and morphisms can only have even colors.

We let $\mathcal{C}',$  and $\mathcal{C}_p'$ denote the subcategories of these categories with the same objects, but only targeted morphisms.

 \section{a variant of the Blanchet-Habegger-Masbaum-Vogel TQFTs}

Blanchet, Habegger, Masbaum, and Vogel   \cite{BHMV2} constructed
3-dimensional TQFT's $(Z_p,V_p)$ on a 3-dimensional cobordism
category where instead of using Lagrangians and integer weights,
as above, to remove anomalies, they used `$p_1$-structures.' 
 These
are versions of the Reshetikhin-Turaev TQFT associated to $SU(2)$
(when  $p$ is even) and to SO(3) (when $p$ is odd.) 

 One can construct by the same  procedure TQFTs (in the sense of
\cite{BHMV2}) defined on ${\mathcal{C}_{\breve p}},$ if
$p$ is odd and on ${\mathcal{C}_{p}},$ when $p$ is even. One can also
construct a quantization functor on ${\mathcal{C}_{p}},$ if $p$ is
odd which satisfies all the  \cite{BHMV2}) axioms except the tensor product
axiom. This failure of the tensor product axiom leads to failure
of the trace property.  We will  denote these related functors on ${\mathcal{C}_{p}}$ also by $(Z_p,V_p).$ We only consider the case $p \ge 3.$ 

These are functors to the category $\FMod(K_p)$ where 
\[
K_p=
\begin{cases}
\BZ [A_p,\frac 1 p],
&\text{if $p \equiv -1 \pmod{4}$ or if $p \equiv 0 \pmod{4}$ }
\\
\BZ [\alpha_p,\frac 1 p],
&\text{if $p \equiv 1 \pmod{4}$ or if $p \equiv 2 \pmod{4}.$ }
\end{cases}
\]
Here $A_p$ is a primitive $2p$th root of unity, and $\alpha_p $ is a primitive $4p$th root of unity, with $A_p= \alpha_p^2.$ Actually $K_3$ could have  been  defined so as to not  include $\frac 1 3$.  We let an element $\kappa_p$ play the role of $\kappa^3$ in \cite{BHMV2}. i.e. raising $w(N)$ by one multiplies $Z_p(N)$ by $\kappa_p.$ We need $4p$th roots of unity in the case  $p$ is $1$ or 2 modulo 4, since  $\kappa_p^2 = A^{-6-\frac{p(p+1)}{2}}.$ 
A basis for $V_p(\Sigma)$ 
is given by  the $p$-colored banded trivalent graphs described in
\cite[4.11,4.14]{BHMV2}. We have bilinear forms (whose adjoints are isomorphisms)

$$v_{\Sigma,p}:V_p(\Sigma ) \times V_p(-\Sigma) \rightarrow K_p,$$
 given by
$$v_{\Sigma,p}([N_{1}],[N_{2}])= <N_{1} \cup_{\Sigma}N_{2}>_p.$$

Let $\eta_p \in K_p$ be the element specified
 in \cite[p.897]{BHMV2} substituting $\kappa_p$ for  $\kappa^3$. 
Let $\Omega_p$ be the linear combination over $\BZ [A_p]$ of $p$-colored cores of a solid torus specified
in \cite[5.8]{BHMV1, BHMV2}.
If $M$ is a closed connected 3-manifold in ${\mathcal{C}_p},$
the associated quantum invariant $<M>_p$
can be described as follows. Let
$L$ be a framed link in $S^3$ that is a surgery description of $M,$ and has signature $w(M),$ and say $n$ components. The $p$-colored banded trivalent graph $G$ can be described by a $p$-colored banded trivalent graph $J$ in the complement of $L.$ Let $L(\Omega_p)$ be the linear combination over $\BZ [A_p]$ that we obtain if we replace each component by $\Omega_p.$
We have that $<M>_p$ is $\eta_p^{1+n}$ times the Kauffman bracket
 of the following $\BZ [A_p]$-linear combination of  the $p$-colored links: the multilinear expansion of $L(\Omega_p) \cup J$ evaluated at $A_p.$
One can see from the long exact sequence of the pair consisting of the 4-ball with 2-handles attached along $L$ relative $M$ that:
\begin{equation}
 n \equiv w(M)+ \be_1(M) \pmod{2}.  \label{closedweight}  \end{equation}
 
 \section{ The  $(S_p,Z_p)$  functors}

{\em In this section, we suppose that  $p$ is either an odd prime or twice an odd prime.}
Consider the rings of integers
\[
 {\De}_p=
\begin{cases}
\BZ [A_p],
&\text{if $p \equiv -1 \pmod{4}$}
\\
\BZ [\alpha_p],
&\text{if $p \equiv 1 \pmod{4}$ or if $p \equiv 2 \pmod{4}.$ }
\end{cases}
\]
These are, of course, Dedekind domains.
Let 
 $\mathcal{D}_p = \eta_p^{-1}$.  $\mathcal{D}_p$ lies in $\kappa_p \BZ [A_p] \subset   {\De}_p.$
 $\mathcal{D}_p$ is a square root of $\frac{-p}{(A_p^2-A_p^{-2})^2}.$
Generalizing earlier work of Murakami, 
 Masbaum and Roberts \cite{MR} give a
nice proof of an integrality theorem which implies that $\mathcal{D}_p <M>_p\in  {\De}_p,$ if $M$ is a closed weighted connected 3-manifold with  a $p$-colored
link. We note that one can further generalize this result to
include $p$-colored banded trivalent graphs.  This is because a $p$-colored banded
graph can be  expanded as a $\BZ[A_p]$-linear combination of
links. This uses the fact that the Temperley-Lieb idempotents can
be defined over $\BZ[A_p]$ \cite{MR}.  It follows that
$(V_p,Z_p)$ is almost integral on $\mathcal{C}_p.$ Thus Theorem
\ref{main} applies.

If $\Si$ is a nonempty object in  ${\mathcal{C}_p}'$,  let $S_p(\Si)$
denote the $ {\De}_p$-submodule of $V_p(\Si)$ generated by
${\mathcal{C}_p}'$-vacuum  states.
Define $S(\emptyset)$ to be  $ {\De}_p$.  According to  Theorem \ref{main}, $S_p(\Si)$ is a finitely generated projective $ {\De}_p$-module. Also $(S_p,Z_p)$ is a functor from ${\mathcal{C}'_p}$ to $\Mod( {\De}_p).$

 \section{ Free modules }

{\em In this section, we suppose that  $p$ is an odd prime.}

We now wish to show that the modules $S_p(\Si)$ are free. In this section we  only succeed in the case that $p= -1 \pmod{4}.$ The case $p= 1 \pmod{4}$ requires material developed in the rest of the paper.

Recall that the class group $C( {\De})$ of a Dedekind domain $ {\De}$ is the quotient
of the group of fractional ideals $F( {\De})$ by the subgroup of principal
fractional ideals $P( {\De}).$ If $j:  {\De} \hookrightarrow  {\De}'$ is an inclusion of Dedekind Domains,
the map $I \mapsto I {\De}'$ induces a map $j_*:C( {\De}) \ra C( {\De}').$
If $\mathcal{M}$ is a finitely generated $ {\De}$-module, there is an invariant $\mathcal{I}(\mathcal{M}) \in C( {\De})$ with the property that
$\mathcal{I}(\mathcal{M})$ is is the identity element of the group if and only if $\mathcal{M}$ is free \cite[Theorem 10.14]{Ja}. Also, $j_*(\mathcal{I}(\mathcal{M}))= \mathcal{I}(\mathcal{M} \otimes_ {\De}  {\De}').$
We thank Jorge Morales for suggesting the following lemma which is familiar to number theorists.
\begin{lemma} \label{inj} If $p$ is an odd prime, the map $C(\BZ [A_p]) \ra C(\BZ [A_p, \frac 1 p])$ is an isomorphism.
\end{lemma}
\begin{proof} Recall that $(p)=(1+A_p)^{p-1}$ and that $(1+A_p)$ is a prime ideal.
 Using  \cite[p. 180 4F]{Rib}, the map $g$ (below) is surjective and its kernel is the infinite cyclic subgroup generated by $(1+A_p).$ It is easy to check that $f$ (below) is surjective.
 The rest of the commutative diagram  with exact rows and columns can be completed by diagram chasing. Here $1$ denotes the trivial multiplicative group.
\[
\begin{CD}
\quad   @.              {1}   @.              {1} @.                    {1} @.                          \quad   \\
@.                       @VVV                       @VVV                    @VVV                            @.     \\
{1} @>>>               <(1+A_p)> @>>>                P(\BZ [A_p]) @>f>>       P(\BZ [A_p, \frac 1 p]) @>>>    {1}   \\
@.                      @VVV                     @VVV                      @VVV                            @.  \\
{1} @>>>          <(1+A_p)> @>>>                F(\BZ [A_p]) @>g>>       F(\BZ [A_p, \frac 1 p]) @>>>    {1}  \\
@.                      @VVV                      @VVV                    @VVV                            @.  \\
\quad @.                {1}    @>>>             C(\BZ [A_p]) @>>>       C(\BZ [A_p, \frac 1 p]) @>>>    {1}   \\
@.                      @.                      @VVV                    @VVV                            @.        \\
\quad.   @.             \quad  @.               {1} @.                    {1}
\end{CD}
\]
\end{proof}
\begin{thm}\label{p-1} If $p$ is an odd prime, and $p \equiv -1 \pmod{4},$ then $S_p(\Si)$ is a free $Z[A_p]$-module.
\end{thm}

\begin{proof} If $p \equiv -1 \pmod{4},$ then $S_p(\Si)\otimes_{Z[A_p]}Z[A_p,\frac 1 p] $ is isomorphic to $V_p(\Si)$
which is free. By Lemma \ref{inj} and the preceding discussion, then  $S_p(\Si)$ is a free $Z[A_p]$-module.

\end{proof}

 \section{An `index two' cobordism subcategory}

As a step towards obtaining a similar result to Theorem \ref{p-1} in the case $p \equiv 1 \pmod{4},$ we introduce a new subcategory of the cobordism category $\mathcal {C}$.

\begin{Def}\label{evendef}
A morphism $N:\Sigma \rightarrow \Sigma'$ in $\mathcal{C}$
 is called even if and only if
$$w(N)\equiv  \dim \left( \jmath_*\el (\Si) + \jmath'_*\el(\Si ') \right) +\beta_1(N)+\beta_0(N)+  \beta_0(\Sigma)+ \frac {\beta_1(\Sigma')} 2
+ \epsilon(N) \pmod{2}$$
where $\jmath:\Sigma \ra N$ and $\jmath':\Sigma' \ra N$  denote the inclusions and $\epsilon (N)$ is one if exactly one of $\Si$  and  $\Si '$ is nonempty. Otherwise 
$\epsilon (N)$ is zero.

\end{Def}

\begin{thm}
 \label{even}
 The composition of two even morphisms is even. 
\end{thm}

Suppose $H$ is a handlebody $H,$  and $\partial H$ is equipped with the Lagrangian given by the kernel of the map induced on the first homology by the inclusion of $\partial H$ into $H.$ If $H$ is  viewed as a morphism from the empty set to $\partial H$ and $H$ is given the weight zero, then $H$ is an even morphism.
This motivated the definition of even morphism.
 We note that if we deleted the $\epsilon(N)$ term in the above definition, then
Theorem \ref{even} would still be true. However the morphism $H$ discussed above
would no longer be even. The proof of Theorem  \ref{even} as given below requires the use of rational Lagrangian subspaces rather than real Lagrangian subspaces. This is why we departed from Turaev and Walker in using rational Lagrangian subspaces for extra structure. Recently,  together with   Qazaqzeh \cite{QG}, we  have obtained the conclusion of Theorem \ref{even}  for cobordisms in the  cobordism  category where surfaces are equipped with real Lagrangians. The  proof relies on Theorem \ref{even} as stated above.

Before we begin the proof, we  reformulate the notion of evenness.
First we specify a homologically standard $\Si$-coboundary.
Given a nonempty extended surface $\Si,$
build $P(\Si)$
(non-uniquely) in four steps:
\begin{enumerate}
\item Form $\Sigma \times I.$
\item Add 1-handles along $\Sigma \times {1}$ if necessary so
that
$\partial_+(\text {result of step 1})$ is connected.
\item Add 2-handles to the result of step 2 along a
collection of disjoint simple curves  in  $\partial_+(\text{ result
of step 2})$ representing some basis for $\el(\Si).$
\item Add 3-handles to the result of step 3 along the
2-sphere components of $\partial_+(\text{ result of step 3}).$
\end{enumerate}

\noindent We use $\Si$ to denote the subset  $\Si \times{0}$ of  $P(\Si)$. This should cause no confusion.
Let $i_{\Sigma}:
 \Sigma \rightarrow P(\Sigma)$ be the inclusion.
Note that $H_0(P(\Sigma)) \approx \mathbb{ Z},$ and
$H_2(P(\Sigma)) \approx \mathbb{ Z}^{\beta_0({ \Sigma})-1}.$
We have
$H_1(P(\Sigma))$ is free Abelian with rank ${\frac 1 2 \beta_1(\Sigma)},$ and $\el(\Si)$ is the kernel
of $i_{\Sigma *}:H_1(\Sigma) \rightarrow H_1(P(\Sigma)).$ If $\Si = \emptyset$, define
$P(\Si) = \emptyset$ as well.
If $N:\Sigma \rightarrow \Sigma'$ is  a morphism in $\mathcal{C},$ let
$$K(N) = -P(\Si) \cup_{\Si} N  \cup_{\Si'}P(\Si').$$ We think of $K(N)$ as $N$ capped off.
Note that the rational homology of $K(N)$ is independent of the choices made.
In the proofs of Proposition \ref{evenprop} and Theorem \ref{even} below,  all homology groups will be taken with rational coefficients.

\begin{prop} \label{evenprop}
$N:\Sigma \rightarrow \Sigma'$ is  even if and only if
$$w(N)\equiv \beta_0(\Sigma')+ \frac {\beta_1(\Sigma')} 2
 + \beta_0(K(N) ) + \beta_1(K(N) ) \pmod{2}.$$
\end{prop}

\begin{proof} 
Let $P$ denote $P(\Si) \sqcup P(\Si').$ Thus $\beta_0(P)=1$ if and only if 
 $\epsilon (N)$ is one. Otherwise $\beta_0(P)$ is even.  Note that $N \cap P = \Si \sqcup \Si'$, so we can identify
$H_i(N \cap P)$ with $H_i( \Si ) \oplus  H_i( \Si' )$. Consider the Mayer-Vietoris sequence for 
$K(N)$ as the union of $N$ and $P.$
\begin{multline*}  H_1( \Si ) \oplus  H_1( \Si' ) \ra 
H_1(N )\oplus H_1\left( P(\Si)  \right) \oplus  H_1\left( P(\Si')  \right) 
 \ra H_1(K(N) )\\
 \ra H_0( \Si ) \oplus  H_0( \Si' )\ra H_0( N  ) \oplus  H_0( P ) \ra H_0(K(N) ) \ra 0
 \end{multline*}
Let $\delta$ denote the dimension of the cokernel of the first map, which is given by
 $ \beta_1(N)-\dim \left( \jmath_*\el (\Si) + \jmath'_*\el(\Si ') \right).$ We have that
$$ \delta + \beta_0(N) + \beta_0(\Sigma) + \epsilon(N) \equiv
\beta_0(\Si')  +  \beta_0(K(N)) + \beta_1(K(N)  \pmod{2}$$
The result now follows.

\end{proof}

\begin{proof} [ Proof of Theorem \ref {even}]  Let $N_1:\Si \ra \Si',$ and $N_2:\Si' \ra \Si''$ be even morphisms in
$\mathcal{C}.$ We wish to show that $N=N_2 \circ N_1:\Si \ra \Si''$ is even. We may reduce to the case that $\be_0(\Si')\ne \emptyset.$  
Let $W_i$ denote a simply connected 4-manifold with boundary $K(N_i),$ and signature $w(N_i).$  Let $W= W_1 \cup_{P(\Si')} W_2.$
From the long exact sequences of the pairs for $(W_i,K(N _i))$ and $(W,K(N) ),$
and from the 4-manifold interpretation of the weight of a composition, we have that
 $$\be_2(W_i)  \equiv w(N_i)+ \be_1(K(N_i ) ) \pmod{2} \  \text{and}
\  \be_2(W) \equiv w(N)+ \be_1(K(N) )  \pmod{2}$$

Since each $N_i$ is even,  by Proposition \ref {evenprop}, we have that
\begin{equation} \label{even2}\be_2(W_1) \equiv \be_0(\Si')+ \frac {\be_1(\Si')}2+ \be_0(K(N_1) ) \pmod{2}, \text{ and }
\end{equation}
\begin{equation} \label{even3} \be_2(W_2) \equiv \be_0(\Si'')+ \frac {\be_1(\Si'')}2+\be_0(K(N_2) ) \pmod{2}
\end{equation}
To show that $N$ is even,  it suffices, by Proposition \ref {evenprop},  to show that
\begin{equation} \label{even4}
\be_2(W) \equiv \be_0(\Si'')+ \frac {\be_1(\Si'')}2+ \be_0(K(N) ) \pmod{2}.
\end{equation}

Because $P(\Si')$ is connected and each component of $W_1$ and  $W_2$ is simply connected, we conclude using a Mayer-Vietoris sequence that $\be_1(W)=0.$ Thus by Lefschetz
 duality and the universal coefficients theorem, 
$$\be_3(W,K(N) )=0 \text{ and } \be_3(W_i,K(N_i) )=0.$$

Therefore, using the long exact sequence for the pairs $(W,K(N)), $ $(W_i,K(N_i)), $
\begin{equation} \label{even5}
\be_3(W)= \be_0(K(N) )-1  \text{ and } \be_3(W_i)= \be_0(K(N_i) )-1 
\end{equation}

We have the Mayer-Vietoris sequence for $W$ as the union of $W_1$ and $W_2$:
\begin{multline} \label{even7} 0\ra H_3(W_1 ) \oplus H_3(W_2 ) \ra H_3(W ) \\
 \ra H_2(P(\Si') )\ra H_2(W_1 )\oplus H_2(W_2 )\ra H_2(W) \ra H_1(P(\Si') )
\ra 0   \end{multline}
Thus by  the exactness of sequence \eqref{even7}
\begin{equation} \label{even8}
 \be_2(W)\equiv \be_2(W_1)+\be_2(W_2)+ \be_2(P(\Si'))+\be_1(P(\Si')) +\be_3(W_1) + \be_3(W_2) +\be_3(W) \pmod{2}\end{equation}

Recall we also have that
 \begin{equation} \label{even9}\be_2(P(\Si'))=  \be_0(\Si')-1 \quad
\text{and} \quad
\be_1(P(\Si'))= \frac {\be_1(\Si')}2.\end{equation}
\noindent
We can derive equation \eqref{even4}, by substituting  equations
\eqref{even2},
\eqref{even3},
\eqref{even5},
 and
\eqref{even9}
into equation \eqref{even8}.

\end{proof}

Note that the identity, $\text{id}_{\Si}: \Si \ra \Si,$ is an even morphism.
This follows from Theorem \ref{even}, but also directly from the definition.
Note the asymmetric treatment
of $\Sigma$ and $\Sigma'$  in Definition \ref{evendef}. Thus the oriented manifold of $N$ viewed as
a morphism from
$-\Sigma'$ to $-\Sigma$ may have a different parity than it does as
originally viewed
as a morphism from $\Sigma$ and $\Sigma'.$ In other words the choice
of incoming and outgoing manifold is important data in a cobordism, which
may affect its parity.

If $N: \Si \ra \Si,$ let $\check{N}:\emptyset \ra -\Si \sqcup \Si$ denote the same manifold with the same weight
and the same Lagrangian but with a different choice of source and target. Similarly define
$\hat{N}: -\Si \sqcup \Si   \ra \emptyset.$

\begin{prop}  \label{evenslice}   The closure  $C(N)$ is the composition $ \widehat{\text{id}_{\Si}} \circ \check{N}.$\quad
 $N: \Si \ra \Si$ is even if and only if $C(N)$ is even.
\end{prop}
\begin{proof}
The first statement follows from the 4-manifold interpretation of
weights. The second statement holds because the parity change
between $N$ and $\check{N}$ is the same as the parity change
between ${\text{id}_{\Si}}$ and $\widehat{\text{id}_{\Si}}.$
\end{proof}

As a corollary,  we  obtain the following congruence which is needed in \cite{GKP}.

\begin{cor} \label{defectmod2} If $ \tilde M \rightarrow M$ is a simple $Z_{d}$-covering of connected 3-manifolds, then
$$\text{ def }( \tilde M \rightarrow M)\equiv  d\left(\beta_1( M)+1\right)+ \beta_1( \tilde M)+1 \pmod{2}.$$\end{cor}
\begin{proof}  It suffices to prove this when $M$ is connected. If we equip $M$ with the weight $\beta_1(M)+1,$ then $M$ is even.
So it is the closure of an even morphism $N.$ So $\tilde M$ is the closure
of the even morphism $N^d.$ It follows that $w(\tilde M)=\beta_1( \tilde M)+1 \mod (2).$ The result then follows from Proposition \ref{defect}. \end{proof}

We let $\mathcal{C}_+ $   and ${\mathcal{C}_p}_+$ 
denote the subcategories of  $\mathcal{C}$   and ${\mathcal{C}_p}$, with the same objects but only  even morphisms.  
We let $\mathcal{C}_+'$  and ${\mathcal{C}_p}_+'$ denote the subcategories of $\mathcal{C}_+$ and  ${\mathcal{C}_p}_+$
with the same objects, but only targeted morphisms.
In other words: $+$ means `take only even morphisms,'  prime means `only use targeted morphisms,'  and  $p$ means `decorate with $p$-colored  banded trivalent graphs.' 

If one restricts
to this `index two' cobordism subcategory, one 
obtains functors, related to the TQFT functors defined by Turaev
with initial data a modular category, but without taking a
quadratic extension of the ground ring of the modular category as
is sometimes needed in \cite[p.76]{Tu}.

\section{  The  $(V^+_p,Z^+_p)$ functors}

{\em In  this section we assume that $p \equiv 1 \pmod{4}$ or $p \equiv 2 \pmod{4}.$}

\begin{prop}  \label{closedeven} If $M$ is an even closed 3-manifold in ${\mathcal{C}_p},$  then $<M>_p \in \BZ [A_p,\frac 1 p].$
If $M$ is an odd closed 3-manifold in ${\mathcal{C}_p},$  then $<M>_p \in \alpha_p \BZ [A_p,\frac 1 p]. $
\end{prop}

\begin{proof}
It is enough to prove this for $M$ connected. If $M$ is even, then $w(N) \equiv 1+ \be_1(M) \pmod{2}.$ By equation \eqref{closedweight}, $n$ is odd. So $\eta^{1+n}\in \BZ [A_p,\frac 1 p].$ The second statement is proved similarly.
\end{proof}

If $\Si$ is an object of ${\mathcal{C}_p}_+,$ let $V_p^+(\Si)$ denote the $\BZ [A_p,\frac 1 p]$ submodule of $V_p(\Si)$ generated by ${\mathcal{C}_p}_+$-vacuum states. Theorem \ref{even} implies:

\begin{prop}
If $N:\Si \ra \Si'$ is a morphism of ${\mathcal{C}_p}_+,$ then $Z_p(N)(V_p^+(\Si))
\subseteq V_p^+(\Si').$
\end{prop}

\begin{prop}
$V_p^+(\Si)$ is a free $\BZ [A_p,\frac 1 p]$-module.
Also $V_p^+(\Si)$ has the same rank as
$V_p(\Si.)$
\end{prop}
\begin{proof}
 The pairing
$v_{\Si,p}: V_p(\Sigma) \otimes V_p(-\Sigma) \ra \BZ [\alpha_p,\frac 1 p]$
has an adjoint:
$$a_{\Si,p}:V_p(\Sigma) \ra \hom\left(V_p(-\Sigma) , \BZ [\alpha_p,\frac 1 p]\right).$$

Take bases for $V_p(\Sigma)$ and $V_p(-\Sigma)$ consisting of
linearly combinations over $\BZ[A_p,\frac 1 p]$ of  ${\mathcal{C}_p}_+$-vacuum states.
Consider the matrix for $a_{\Si,p}$
with respect to these bases. By Proposition \ref{closedeven}, the entries of the matrix all lie in $\BZ [A_p,\frac 1 p].$  Let $\gamma$  denote the determinant of this matrix. $\gamma$ is a unit from $\BZ [\alpha_p,\frac 1 p],$ and
 $ \gamma \in \BZ [A_p,\frac 1 p].$ As  $\gamma^{-1} \in \BZ [\alpha_p,\frac 1 p],$ we conclude
that $\gamma^{-1} \in \BZ [A_p,\frac 1 p].$ In other words, $\gamma$ is a unit in $\BZ [A_p,\frac 1 p].$

Let $W(\Sigma)$  and $W(-\Sigma)$ denote the $\BZ [A_p,\frac 1 p]$ submodules spanned by the  chosen bases for $V_p(\Sigma)$ and $V_p(-\Sigma).$ We have that $V^+_p(\Sigma) \supseteq W(\Sigma).$  The restriction of $a_{\Si,p}$ to
 $V^+_p(\Sigma)$ is an injective map $V^+_p(\Sigma) \ra \hom\left(W(-\Sigma) , \BZ [A_p,\frac 1 p]\right).$  But its further restriction to
$W(\Si)$ is an isomorphism as it is given by a matrix over $\BZ [A_p,\frac 1 p]$ whose determinant is a unit. It follows that $V^+_p(\Sigma)$ equals $W(\Si)$ and so $V^+_p(\Sigma)$ is free.
 \end{proof}

Thus we obtain a functor $(V_p^+,Z_p)$ from $\mathcal{C}^+_p$ to $\FMod(\BZ [A_p,\frac 1 p])$.

\section{ The  $(S^+_p,Z^+_p)$ functors}

{\em In this section, we suppose that  $p$ is either an odd prime  and $p \equiv 1 \pmod{4}$ or $p$ is twice an odd prime.}

Let 
$\mathcal{D}^+_p= \kappa_p \mathcal{D}_p  \in \BZ [A_p].$
For a closed connected  3-manifold $M$ of $\mathcal{C}_{p+},$
$\mathcal{D}^+_p <M>_p \in \BZ [A_p].$
So
$(V_p^+,Z_p)$ is almost
$\BZ [A_p]$-integral.
If $\Si$ is a nonempty object in  $\mathcal{C}_{p+}',$ let $S^+_p(\Si)$
denote the $\BZ [A_p]$-submodule of $V_p(\Si)$ generated by $\mathcal{C}_{p+}'$-vacuum states. Define $S^+_p(\emptyset)$ to be the $\BZ [A_p]$-submodule of $V_p^+(\emptyset)=K_p$ generated by $1.$ By   Theorem \ref{main}, $S^+_p(\Si)$ is a finitely generated projective $\BZ [A_p]$-module. Also $(S^+_p,Z_p)$ is a functor from $\mathcal{C}'_{p+}$ to $\Mod(\BZ [A_p]).$
When $p$ is twice an odd prime, we are unable, at this point, to prove  $S^+_p(\Si)$ is a free. 

\begin{thm} 
If $p$ is an odd prime, and $p \equiv 1 \pmod{4},$ then  $S^+_p(\Si)$ is a free $Z[A_p]$-module. It follows that $S_p(\Si)$ is also a free $Z[\alpha_p]$-module.

\end{thm}
\begin{proof} 
If $p \equiv 1 \pmod{4},$ $S^+_p(\Si)\otimes_{Z[A_p]}Z[A_p,\frac 1 p] $ is isomorphic to $V^+_p(\Si)$
which is free. By Lemma \ref{inj} and the preceding discussion, $S^+_p(\Si)$ is a free $Z[A_p]$-module.
$S_p(\Si)$ is isomorphic to $S^+_p(\Si)\otimes_{Z[A_p]}Z[\alpha_p],$ and so is a free $Z[\alpha_p]$-module.
\end{proof}

\section{ Strong Shift equivalence Class Invariants}

Strong shift equivalence is a concept due to R.F. Williams which comes from symbolic dynamics \cite{LM,Wag}. Originally defined for square matrices with non-negative entries, this notion is meaningful in any category. An elementary equivalence between two morphisms $E_1$ and $E_2$ in a category is a pair of morphisms $(R,S)$ in the category such that $E_1 = RS,$ and $E_2= SR.$
Strong shift equivalence is the relation on morphisms in the category that is generated by elementary equivalences.

$\mathcal{C}$ and $\mathcal{C^+}$ have the connected slitting property.
The decorated categories also have this  property.
We  state a version of a previous theorem of ours from \cite{SSE} in this context:

\begin{thm}  \label{infinite} Let $M$ be  a closed connected morphism in $\mathcal{C}$,  and $\chi:H_1(M) \rightarrow \mathbb{ Z}$ be an epimorphism.
Any two connected fundamental domains for the infinite cyclic cover of $M$ specified by $\chi$ are strong shift equivalent in $\mathcal{C}'$. 
Suppose $(V,Z)$ is almost $ {\De}$-integral on $\mathcal{C}$.  If $N$ is a connected fundamental domain for the infinite cyclic cover of $M$ specified by $\chi,$ then
 $Z(N):S(\Si)\ra S(\Si)$ is well defined up to strong shift equivalence in $\Mod( {\De}).$
\end{thm}

There are also versions of this theorem for decorated cobordism categories.  We also have a version  for the even cobordism category.

Knot invariants result when we take $M$ to be 0-framed surgery along the knot, and  $\chi$ to be the homomorphism which sends a positive meridian to one. If one then applies the functor $(S_p,Z_p)$, one obtains knot invariants taking values in the strong shift equivalence classes of  morphisms in  $\FMod(D_p).$    This was one of the motivations for developing $(S_p,Z_p).$ In particular one may now define the analog of the Bowen-Franks invariant in this situation \cite{LM}. The Bowen-Franks invariant of a knot is the  $ {\De}$-module given by $\text{the cokernel of } Z(N)-1:S(\Sigma) \rightarrow S(\Sigma) , $ where $N$ is a fundamental  domain for the infinite cyclic cover of  0-framed surgery along the knot.  We  hope to investigate this invariant more elsewhere.

\section{ simple covers using $(V_p,Z_p)$}

{\em In this section $p$ is an odd prime.}

If $M$ is a closed 3-manifold, with a $p$-colored banded trivalent graph $G$,  let
$$\langle\langle M \rangle  \rangle_p=  \langle \flat M \rangle _p,$$
where $\flat M$ is $M$ with integral weight zero.  If $M$ already has a weight $w(M),$ we have
$\langle\langle   M \rangle  \rangle _p= \kappa_p^{-w(M)} \langle M \rangle _p.$ Note the graph $G$ is contained in the data of $M.$
 See \cite{per} for some results related to the first part of the following theorem when $r \ne p.$

\begin{thm}  \label{simplecong}
If $M$ is a simple $Z_{d}$  covering space of $M$ with a $p$-colored banded trivalent graph $G$ where $d=r^s,$ and  $r$ is prime,
then
$$ \langle\langle\tilde M \rangle  \rangle _p - \kappa_p^{\text{ def }(\tilde M \rightarrow M)}\left(\langle\langle M \rangle  \rangle _p\right)^d \in r {\De}_p.$$
In particular, if  $p=r,$ then for some  $n \in \BZ,$ we have:
$$ \langle\langle\tilde M \rangle  \rangle _p - 
i^{ \frac {p+1} 2\left( \beta_1(M) +\beta_0(M)\right)}
\kappa_p^{\text{ def }(\tilde M \rightarrow M)} n \in p {\De}_p.$$
\end{thm}

\begin{proof} First we consider the case when all the colors of $G$ are even.
Then we can work in a situation where $(V,Z)$ has the trace property.
Assign $M$ the weight $ \beta_1(M) +\beta_0(M),$ so that $M$ is  even. Applying  Corollary \ref{scover} with $ {\De}=  {\De}_p$, and $(V,Z)= (V_p,Z_p),$
we obtain
$ \langle\tilde M \rangle _p-\left(\langle M \rangle _p\right)^d \in r {\De}_p,$ or
$$\kappa_p^{w(\tilde M)} \langle\langle\tilde M \rangle  \rangle _p-\left( \kappa_p^{w(M)} \langle\langle M \rangle  \rangle _p \right)^d \in r {\De}_p.$$ By Proposition \ref{defect}, $dw( M)- w(\tilde M)=\text{ def }(\tilde M \rightarrow M).$  Now suppose that $p=r$.
If $M$ is even,  then $ \langle M \rangle  \in \BZ[A_p].$ Thus
$ \langle \langle M \rangle  \rangle _p \in  \kappa_p^{ w(M)}\BZ[A_p]=
 i^{ \frac {p+1} 2 w(M)}\BZ[A_p].$
If  $d=p^s,$ then $\left( \langle  \langle M \rangle  \rangle _p\right)^d$ differs from  an integer times $i^{ \frac {p+1} 2 w(M)}$ by an element in $r {\De}_p.$

Now suppose some colors of $G$ are odd. We can replace $G$ by a
$\BZ[A_p]$-linear combination  of banded links (colored one)
without changing $ \langle  \langle M \rangle  \rangle _p.$ Next replace $\tilde G$ by same
linear combination of the inverse image of these links. This will
not change $ \langle  \langle \tilde M \rangle  \rangle _p.$ This makes use of the
idempotent property of the Temperley-Lieb idempotents. Next we
recolor each link component $p-3$ without changing the quantum
invariants upstairs or downstairs \cite[(6.3)iii]{BHMV1}. The
result now follows by the earlier case.
\end{proof}

\section{ A conjecture}

We close this paper with a conjecture that the above work suggested. Let  $$I_p(M)= \mathcal{D}_p  \langle  \langle M \rangle  \rangle _p \in  {\De}_p.$$

The inner rank of a group  is the maximal $k$ such that there is an epimorphism of that group onto the free group on $k$ letters \cite{Ly}.
If $M$ is an oriented connected
three manifold, define $c(M),$ the cut number of $M,$ to be the
maximal number of oriented surfaces one can place in $M$ and still
have a connected complement.
Jaco \cite{Jaco} showed that the inner rank of $\pi_1(M)$ is the cut number of $M.$
See also \cite[Theorem 7.4.4]{Kai} for a discussion of these matters.

\begin{con} Let $M$ be an oriented connected closed three manifold (with possibly a banded trivalent colored graph), and $p$ be a prime,
then ${\mathcal{D}_p}^{c(M)}$ divides $I_p(M).$
\end{con}

The conjecture for $p=5$ was obtained (in case that the graph is empty) together
with   Kerler using Kerler's basis. Just recently,  this conjecture has been obtained together with  Masbaum  for all odd primes  $p$ (again assuming the graph is empty) \cite{GM}.
Cochran and Melvin \cite{CM} have studied the divisibility of $I_p(M)$ by $\mathcal{D}_p.$
According to  \cite{CM} if   $M$ is zero framed surgery along a boundary link with $c$ components then  ${\mathcal{D}_p}^{c}$ divides $I_p(M).$ The conjecture also implies this result.

\end{document}